\documentclass{article}

\usepackage{amsfonts}
\usepackage{amssymb}
\usepackage{amsmath}
\usepackage{amsthm}
\usepackage{cite}

\newcommand{\C}{{\mathbb C}}
\newcommand{\N}{{\mathbb N}}

\newcommand{\lcm}{{\operatorname{lcm}}}
\newcommand{\tp}{^{\rm t}}
\newcommand{\rd}{{\mathfrak r}}

\newcommand{\rg}{{r}}

\newcommand{\bcdot}{$\discretionary{\mbox{$ \cdot $}}{}{}$}

\newtheorem{stelling}{Theorem}[section]
\newtheorem{propositie}[stelling]{Proposition}
\newtheorem{lemma}[stelling]{Lemma}

{
\theoremstyle{definition}

}

\newenvironment{bewijs}{\begin{proof}}{\end{proof}}

\title{Another generalization of Mason's ABC-theorem}
\author{Michiel de Bondt}

\begin{document}

\maketitle

The well-known ABC-conjecture is generally formulated as follows:

\paragraph{The ABC-conjecture.}
\begin{em} 
Consider the set $S$ of triples $(A, B, C) \in \N^3$ such that 
$ABC \ne 0$, $\gcd\{A,B,C\} = 1$ and
$$
A + B = C
$$
Then for every $\epsilon > 0$, there exists a constant $K_{\epsilon}$
such that
$$
C \le K_{\epsilon} \cdot R(ABC)^{1+\epsilon}
$$
for all triples $(A, B, C) \in S$, where $R(ABC)$ denotes the square-free
part of the product $ABC$.
\end{em}\par \bigskip

The ABC-conjecture is studied in many papers, and this article will
not be another of them. Instead, we consider an analog of this conjecture 
for polynomials over $\C$ instead of integers: Mason's ABC-theorem:

\paragraph{Mason's ABC-theorem.}
\begin{em} 
Let $f_1, f_2, f_3$ be polynomials 
over $\C$ without a common factor, not all constant, such that 
$$
f_1 + f_2 + f_3 = 0
$$
Then 
$$
\max_{1 \le m \le 3} \deg f_m \le \rg(f_1f_2f_3) - 1
$$
where $\rg(g)$ denotes the number of distinct zeros of $g$.
\end{em}\par\bigskip

This theorem was proved at first by Stothers in \cite{St}. So Mason did what 
Stayman did with the bridge convention that has his name: he made the theorem
known, even popular. 

The bound in Mason's theorem can be reached by examples of arbitrary
large degree, namely $f_1 = f^3, f_2 = {\rm i}g^2, f_3 = -(f^3-g^2)$, where
$f$ and $g$ reach H. Davenport's bound:
$$
\deg (f^3 - g^2) \ge \frac12 \deg f + 1
$$
All $f$ and $g$ that reach the Davenport bound are determined in \cite{ZaD}.
The easiest example is
$$
(x^2 + 2)^3 - (x^3 + 3x)^2 = 3x^2 + 8
$$

So Mason's theorem seems the best you can get. But there is room for
generalization. One direction is followed for the ABC-conjecture as well,
namely adding more integers/polynomials to (get) the sum that vanishes.
Another direction is allowing more indeterminates in the polynomials. We will
discuss both generalizations. There has already been done a lot of work
in these direction, mainly using so called {\em Wronskians}, 
but it seems that no one has combined all ideas to get the best generalized 
results one can get by means of Wronskians.

A third direction of generalization is to use elements of so-called 
{\em function fields} instead of univariate polynomials 
\cite{BrMa,Za,HY1}, or using {\em meromorphic functions} instead of 
multivariate polynomials \cite{HY}. 
These generalizations will decrease the readability of this expository 
paper, so we restrict ourselves to polynomials.

\section{Generalizations of Mason's ABC-theorem} \label{mason34}

Let $p$ be a (possibly multivariate) polynomial over $\C$. 
Then we can factorize $p$: 
$$
p = p_1^{e_1} p_2^{e_2} \cdots p_s^{e_s}
$$
with all $p_i$ irreducible and pairwise relatively prime, 
and all $e_i \ge 1$. Let
$$
\rd(p) := p_1 p_2 \cdots p_s
$$
be the square-free part of $p$ and denote by $\rg(p)$ the degree of $\rd(p)$.

Associating polynomials with principal ideals, we have that $\rd(p)$
is the radical of $p$; hence the symbol $\rd$ is used.

Mason's ABC-theorem for three polynomials is generally formulated as follows
\cite{Mason,Sil,Sny,St}:

\begin{stelling} \label{mason3}
Let $f_1, f_2, f_3$ be pairwise relatively prime univariate polynomials 
(in the same variable) over $\C$, not all constant, such that 
$$
f_1 + f_2 + f_3 = 0
$$
Then 
$$
\max_{1 \le m \le 3} \deg f_m \le \rg(f_1f_2f_3) - 1
$$
\end{stelling}

In \cite[Theorem 1.2]{ShSp}, H.N. Shapiro and G.H. Sparer generalize theorem 
\ref{mason3} as follows, see also \cite{HY}:

\begin{stelling} \label{masons}
Let $n \ge 3$ and $f_1, f_2, \ldots, f_n$ be pairwise relatively prime 
(possibly multivariate) polynomials over $\C$, not all constant, 
such that 
$$
f_1 + f_2 + \cdots + f_n = 0
$$
Then 
$$
\max_{1 \le m \le n} \deg f_m \le (n-2)\Big(\rg(f_1 f_2 \cdots f_n) - 1\Big)
$$
\end{stelling}

In \cite[Theorem 5]{BT}, M. Bayat and H. Teimoori formulate the following
improvement of the estimation bound of theorem \ref{masons} 
(so with all $f_i$s pairwise relatively prime)
as follows: they replace $(n-2)(\rg(f_1 f_2 \cdots f_n) - 1)$ by
$$
(n-2)\left(\rg(f_1 f_2 \cdots f_n) - \frac{n-1}2\right)
$$ 
for the case that at most
one of the $f_i$s is constant and by 
$$
(n-k-1)\left(\rg(f_1 f_2 \cdots f_n) - \frac{n-k}2\right)
$$ 
for the case that exactly
$k \ge 1$ of the $f_i$s are constant. This is indeed
an improvement, for if $k \ge 1$ of the $f_i$s are constant and $k < n$, then 
$n-k-1 \le n-2$ and
$$
\rg(f_1 f_2 \cdots f_n) \ge n-k \ge \frac{n-k}2 \ge 1
$$
because there cannot be exactly one $f_i$ that is not constant

Unfortunately, the proof of \cite[Theorem 5]{BT} is incorrect: 
\cite[Lemma 4]{BT} has counterexamples. But we shall see that the theorem 
itself is correct. In \cite{HY1}, the univariate case of theorem \ref{masons}
is proved, and also the erratic \cite[Theorem 5]{BT} can be viewed as a correct
proof for the univariate case.

But let us first discuss the condition that the $f_i$s are pairwise 
relatively prime. This condition is quite restrictive, so it is a good idea 
to try and get rid of it, and replace it by something
weaker. The example $n=3$, $f_1 = f_2 = x^{100}$, $f_3 = -2 x^{100}$ shows 
that we cannot just forget the condition that all $f_i$s are relatively prime.
So let us replace it by the condition that just
\begin{equation} \label{gcdcond}
\gcd\{f_1, f_2, \ldots, f_n\} = 1
\end{equation}
Now theorem \ref{masons} remains valid for $n=3$, because the conditions
$\gcd\{f_1,\allowbreak f_2,\allowbreak f_3\} = 1$ and 
$f_1 + f_2 + f_3 = 0$ {\em imply} that $f_1, f_2, f_3$ are pairwise 
relatively prime. 

This is no longer the case if $n \ge 4$. Reading the proof of theorem 
\ref{masons} above as given in \cite{ShSp}, it seems that 
$\rg(f_1 f_2 \cdots f_n)$ is just a shorthand
notation for $\rg(f_1) + \rg(f_2) + \cdots + \rg(f_n)$, but if the $f_i$s 
are not pairwise relatively prime, then both expressions are different. So we replace
$\rg(f_1 f_2 \cdots f_n)$ by $\rg(f_1) + \rg(f_2) + \cdots + \rg(f_n)$ as 
well. There are, however, also generalizations with $\rg(f_1 f_2 \cdots f_n)$,
which we will discuss later.

Now the example $n=4$, $f_1 = -f_2 = x^{100}$, $f_3 = -f_4 = (x+1)^{100}$ shows us
that we are not ready yet to prove something. The problem is that
$f_1 + f_2 + \cdots + f_n$ has a proper subsum that vanishes. 
Actually, such proper subsums can be seen as instances of the original sum with
smaller $n$, and it seems reasonable that (\ref{gcdcond}) is satisfied
for these subsums as well, i.e. 
$$
f_{i_1} + f_{i_2} + \cdots + f_{i_s} = 0 \Longrightarrow 
\gcd\{f_{i_1}, f_{i_2}, \ldots, f_{i_s}\} = 1
$$
where $1 \le i_1 < i_2 < \cdots < i_s \le n$. This way we get a valid 
assertion:

\begin{stelling} \label{masons2}
Let $n \ge 3$ and $f_1, f_2, \ldots, f_n$ be (possibly multivariate)
polynomials over $\C$, not all constant, such that 
$$
f_1 + f_2 + \cdots + f_n = 0
$$
Assume furthermore that for all $1 \le i_1 < i_2 < \cdots < i_s \le n$,
$$
f_{i_1} + f_{i_2} + \cdots + f_{i_s} = 0 \Longrightarrow 
\gcd\{f_{i_1}, f_{i_2}, \ldots, f_{i_s}\} = 1
$$
Then 
\begin{equation} \label{kk}
\max_{1 \le m \le n} \deg f_m \le 
   (n-2)\Big(\rg(f_1) + \rg(f_2) + \cdots + \rg(f_n) - 1\Big)
\end{equation}
\end{stelling}

If we replace the constant term $-1$ on the right hand side of (\ref{kk})
by $+n$, then the case in which the $f_i$s are univariate without a vanishing 
proper subsum of $f_1 + f_2 + \cdots + f_n$ follows from 
\cite[Th.\@ B]{BrMa} and the proof of \cite[Cor.\@ II]{BrMa}. An improvement
of the proof of \cite[Cor.\@ II]{BrMa} as indicated in section \ref{masondisc}
below subsequently replaces the term $+n$ by $+(n-1)/2$.

If one does not wish to replace $\rg(f_1 f_2 \cdots f_n)$ by 
$\rg(f_1) + \rg(f_2) + \cdots + \rg(f_n)$ (and neither requires the $f_i$s 
to be prime by pairs), then one can use
the inequality $\rg(f_i) \le \rg(f_1 f_2 \cdots f_n)$ to obtain a coefficient
$n(n-2)$, but in \cite{Vo} and \cite[Cor.\@ I]{BrMa}, it is shown that
in the univariate case, $(n-1)(n-2)/2$ is enough and that 
$-1$ can be maintained within the parentheses. We will prove the multivariate
version of this result:

\begin{stelling} \label{masons2a}
Under the conditions of theorem \ref{masons2},
\begin{equation} \label{kka}
\max_{1 \le m \le n} \deg f_m \le 
   \frac{(n-1)(n-2)}2 \Big(\rg(f_1 f_2 \cdots f_n) - 1\Big)
\end{equation}
\end{stelling}

\section{Improvements of theorems \ref{masons2} and \ref{masons2a}}

But theorems \ref{masons2} and \ref{masons2a} are not the best one can
get. One improvement on \ref{masons2a} is by U. Zannier in \cite{Za}, but
his idea also applies to \ref{masons2}. The
coefficient $n-2$ in (\ref{kk}) should be expressed in the dimension $d$ of the
vector space over $\C$ spanned by the $f_i$s. 
Since $f_1 + f_2 + \cdots + f_n = 0$, $d$ is at most $n-1$, 
so the straightforward improvement is replacing $n-2$ by $d-1$. 
But also the residual term  $(n-2) \cdot -1$ can be improved: 
the natural improvement of the corresponding
term $(n-1)(n-2)/2$ in (9) of \cite[Theorem 5]{BT} is $d(d-1)/2$, so we get 
$$
\max_{1 \le m \le n} \deg f_m \le 
   (d-1)\Big(\rg(f_1) + \rg(f_2) + \cdots + \rg(f_n) - \frac{d}2\Big)
$$

Another improvement is due to P.-C. Hu and C.-C. Yang in \cite{HY1,HY}. 
They extend the definition of the $\rg(g)$ by defining 
$$
\rd_e(g) = \gcd\{g, \rd(g)^e\}
$$
and $\rg_e(g) = \deg \rd_e(g)$. So $\rd_1(g) = \rd(g)$ is the square-free 
part of $g$ and $\rd_2(g)$ is the cube-free part of $g$, etc. Now we have a 
trivial inequality
$$
\rg_e(g) \le e \,\rg(g)
$$
and taking $e = n-2$ indicates precisely how Hu and Yang improve the estimate:
they migrate the coefficient $n-2$ to a subscript of $\rg$. This migration 
has the drawback that the residual term $(n-2) \cdot -1$ does not survive 
several reductions any more (reductions that decrease the dimension of
the vector space over $\C$ spanned by the $f_i$s). 
This can be overcome by only stating that there is a 
$\rho$ with $2 \le \rho \le n-1$, such that
$$
\max_{1 \le m \le n} \deg f_m \le 
   (\rho-1)\Big(\rg(f_1) + \rg(f_2) + \cdots + \rg(f_n) - \frac{\rho}2\Big)
$$
and combining the above idea with that of Zannier, we even assume that 
$\rho \le d$ instead of $\rho \le n-1$.

\begin{stelling} \label{masons3}
Let $n \ge 3$ and $f_1, f_2, \ldots, f_n$ be (possibly multivariate)
polynomials over $\C$, not all constant, such that 
$$
f_1 + f_2 + \cdots + f_n = 0
$$
Assume furthermore that for all $1 \le i_1 < i_2 < \cdots < i_s \le n$,
$$
f_{i_1} + f_{i_2} + \cdots + f_{i_s} = 0 \Longrightarrow 
\gcd\{f_{i_1}, f_{i_2}, \ldots, f_{i_s}\} = 1
$$
Now let $d$ be the dimension of the vector space over $\C$ spanned by the
$f_i$s. Then there exists a $\rho$ with $2 \le \rho \le d$, such that
\begin{align} 
\max_{1 \le m \le n} \deg f_m 
&\le \rg_{\rho-1}(f_1) + \rg_{\rho-1}(f_2) + \cdots + \rg_{\rho-1}(f_n)
      - \frac{\rho(\rho-1)}2 \qquad \label{klhy} \\
&\le (d'-1)\left(\rg(f_1) + \rg(f_2) + \cdots + \rg(f_n) - \frac{d'}2\right)
      \label{kl}
\end{align}
for all $d'$ between $d$ and $n-k+1$ inclusive, 
where $k$ is the number of constant $f_i$s.
\end{stelling}

\begin{bewijs}[Proof of {\cite[Theorem 5]{BT}}]
Since $f_1 + f_2 + \cdots + f_n = 0$, it follows that $d \le n-1$. 
So the first inequality
(9) of \cite[Theorem 5]{BT} follows. Assume that exactly $k$ of
the $f_i$s are constant for some $k$ with $1 \le k \le n-1$, and assume
without loss of generality that $f_n$ is not constant. Since the vector
space over $\C$ spanned by the $k$ constant $f_i$s has dimension $1$ at most,
the vector space over $\C$ spanned by $f_1, f_2, \ldots, f_{n-1}$ has dimension
$(n-1) - (k-1) = n-k$ at most. But since $f_1 + f_2 + \cdots + f_n = 0$, the
latter vector space is also the vector space over $\C$ spanned by 
$f_1, f_2, \ldots, f_n$. So $d \le n-k$ and the second inequality (10) of 
\cite[Theorem 5]{BT} follows as well.
\end{bewijs}

The improvements on theorem \ref{masons2a} are similar to those on theorem
\ref{masons2}:

\begin{stelling} \label{masons3a}
Under the conditions of theorem \ref{masons3}, there exists a $\sigma$
with $1 \le \sigma \le d(d-1)/2$ such that
\begin{align} 
\max_{1 \le m \le n} \deg f_m 
&\le \rg_{\sigma}(f_1 f_2 \cdots f_n) - \sigma \label{klahy} \\
&\le \frac{d'(d'-1)}2 \Big(\rg(f_1 f_2 \cdots f_n) - 1\Big) \label{kla}
\end{align}
for all $d' \ge d$.
\end{stelling}

We postpone the proofs of theorems \ref{masons3} and \ref{masons3a} until
section \ref{masonred}, since we first consider some applications.

\section{Applications to Fermat-Catalan equations}

Just like the ABC-conjecture for integers can be used to tackle 
Fermat's Theorem for integers, versions of Mason's Theorem can be used
to tackle polynomial Diophantic equations:

\begin{stelling}[Generalized Fermat-Catalan] \label{fercat}
Assume
$$
g_1^{e_1} + g_2^{e_2} + \cdots + g_n^{e_n} = 0
$$
and $f_1, f_2, \ldots, f_n$ satisfy the conditions of theorem \ref{masons3}, 
where $f_i = g_i^{e_i}$ for all $i$. Then
$$
\sum_{i=1}^n \frac1{e_i} > \frac1{d-1}
$$
where $d$ is the dimension of the vector space over $\C$ spanned by the
$f_i$s
\end{stelling}

\begin{bewijs}[Proof (based on ideas in \cite{HY1}).]
Assume $f_m$ has the largest degree among the $f_i$s. 
From theorem \ref{masons3}, and $\rg(f_i) \le \deg g_i = e_i^{-1} \deg f_m$, 
it follows that
$$
\deg f_m \le (d-1) \left(\sum_{i=1}^n \frac1{e_i} \deg f_m - \frac{d}2 \right)
$$
which rewrites to
\begin{equation} \label{fchy}
\left(\sum_{i=1}^n \frac1{e_i} - \frac1{d-1}\right) \deg f_m \ge 
\frac{d}2
\end{equation}
which completes the proof.
\end{bewijs}

In \cite[Th.\@ 3.1]{ShSp} and \cite[Th.\@ 8]{BT}, theorem \ref{fercat}
is proved by way of the following inequality:
\begin{equation} \label{fcss}
\left(\sum_{i=1}^n \frac1{e_i} - \frac1{d-1}\right) 
\sum_{i=1}^n \deg g_i \ge \frac{d}2 \sum_{i=1}^n \frac1{e_i}
\end{equation}
but the proof of (\ref{fcss}) will not be copied in a third article today.

In \cite[(3.3)]{ShSp} and \cite[Cor.\@ 10]{BT},
the result of theorem \ref{fercat} is rewritten into a Fermat-type
equation, i.e.\@ with all $e_i$ equal. 
But it is not observed that in the Fermat case,
the condition that the $f_i$s are relatively prime by pairs can be omitted.
Having a version of a generalized Mason's theorem in which the $f_i$s must 
be relatively prime by pairs is only partially an excuse for that, since
it suffices to use the case that $f_1, f_2, \ldots, f_{n-1}$ are
linearly independent of theorem \ref{masons2}, which can be proved
with the methods of \cite{ShSp} and \cite{BT}, 
see also \cite[Th.\@ 1.3]{HY1,HY}.

We say that polynomials $f_1$ and $f_2$ are {\em similar} if
$f_2 = \lambda f_1$ for some $\lambda \in \C^{*}$.

\begin{stelling}[Generalized Fermat] \label{GenFer}
Assume
$$
g_1^e + g_2^e + \cdots + g_n^e = 0
$$
for some polynomials $g_i$, not all zero, and suppose that
$$
e \ge (d+1)(d-1)
$$
where $d$ is the dimension of the vector space over $\C$ spanned by the $g_i^e$s.
Then the vanishing sum $g_1^e + g_2^e + \cdots + g_n^e$ decomposes into 
vanishing subsums
$$
g_{i_1}^e + g_{i_2}^e + \cdots + g_{i_s}^e = 0
$$
with $1 \le i_1 < i_2 < \cdots < i_s \le n$, for which
all $g_{i_j}\!$s are pairwise similar.
\end{stelling}

\begin{bewijs}
Assume without loss of generality that
$$
g_1^e, g_2^e, \ldots, g_d^e
$$
is a basis of the vector space over $\C$ spanned by the $g_i^e$s. It suffices to show that for all $j > d$, $g_j^e$ is similar to $g_i^e$ for some $i \le d$. Assume without loss of generality that
$$
g_j^e = \lambda_1 g_1^e + \lambda_2 g_2^e + \cdots + \lambda_s g_s^e
$$
where $1 \le s \le d$ and $\lambda_1 \lambda_2 \cdots \lambda_s \ne 0$. In order to reduce to the case that the $g_i$s are relatively prime and $d = n-1$, we define
$$
h_i := \frac{\sqrt[e]{\lambda_i} g_i}{\gcd\{g_1, g_2, \ldots, g_s,g_j\}}
$$
for all $i \le s$ and for $i = j$, where $\lambda_j = -1$, since then we get
$$
h_1^e + h_2^e + \cdots + h_s^e + h_j^e = 0
$$
Furthermore, $h_1^e, \ldots, h_s^e$ are linearly independent over $\C$, and
$$
\gcd\{h_1^e,h_2^e,\ldots,h_s^e,h_j^e\} = \gcd\{h_1,h_2,\ldots,h_s,h_j\}^e = 1
$$
If all $h_i$s are constant, then $g_j^e$ is similar to $g_1^e$. So assume that not all $h_i$s are constant. Then it follows from theorem \ref{fercat} that $\frac{s+1}e > \frac1{s-1}$, i.e.\@ $e < (s+1) \bcdot (s-1) \le (d+1)(d-1)$. Contradiction.
\end{bewijs}

\section{A theorem of Davenport} \label{daven}

Now let us look at sums of powers that do {\em not} vanish:
$$
g_1^{e_1} + g_2^{e_2} + \cdots + g_{n-1}^{e_{n-1}} = g_n \ne 0
$$
and suppose that no subsum of $g_1^{e_1} + g_2^{e_2} + \cdots + 
g_{n-1}^{e_{n-1}}$ vanishes. Now the question is how far the degree of $g_n$
can drop. In \cite{Dav}, H. Davenport studied the case
$n = 3$, $e_1 = 3$, $e_2 = 2$, and showed that
$$
\deg(f^3 - g^2) \ge \frac12 \deg f + 1
$$
see also \cite{St}. We shall formulate a generalization of this result that 
improves \cite[(6)]{HY1}, by weakening the conditions.

But first, we need some preparations. In order to get (\ref{kla}) of theorem
\ref{masons3a} from (\ref{klahy}), it suffices to show that
$$
\bigg(\frac{d'(d'-1)}2-\sigma\bigg) \rg(f_1 f_2 \cdots f_n) \ge \frac{d'(d'-1)}2-\sigma
$$
This follows from the fact that not all $f_i$s are constant.
It is somewhat more work to get (\ref{kl}) of theorem \ref{masons3}
from (\ref{klhy}). Since $\rg_{\rho-1}(f_i) \le (\rho-1)\rg(f_i)$ for all $i$, it suffices to show that
$$
\Big((d'-1) - (\rho-1)\Big) \Big(\rg(f_1) + \rg(f_2) + \cdots + \rg(f_n)\Big)
\ge \frac{d'(d'-1)}2 - \frac{\rho(\rho-1)}2
$$
This follows since the right hand side equals $\rho + (\rho+1) + \cdots +
(d'-1) \le ((d'-1) - (\rho-1)) (d' - 1)$ and
$$
\rg(f_1) + \rg(f_2) + \cdots + \rg(f_n) \ge n-k \ge d'-1
$$

Let $k'$ be the number of constant $f_i$s with $i \le n-1$. If $d' \le n - k'$, then
$$
\rg(f_1) + \rg(f_2) + \cdots + \rg(f_{n-1}) \ge n-1-k' \ge d'-1
$$
So for $d'$ with $d \le d' \le n - k'$ (such a $d'$ exists because $d \le n-k \le n-k'$), we do not need $f_n$ to boost the residual term to $d'(d'-1)/2$:
$$
\max_{1 \le m \le n} \deg f_m 
\le (d'-1)\left(\rg(f_1) + \rg(f_2) + \cdots + \rg(f_{n-1}) 
    - \frac{d'}2\right) + \rg_{\rho-1}(f_n)
$$
Estimating $\rg_{\rho-1}(f_n)$ by $\deg f_n$ and realizing that at least
two $f_i$s have maximum degree, we get (\ref{predav}) of theorem 
\ref{predavth} below under the conditions of theorem \ref{masons3}:

\begin{stelling} \label{predavth}
Let $f_1, f_2, \ldots, f_n$ be (possibly multivariate)
polynomials over $\C$, not all constant, such that 
$$
f_1 + f_2 + \cdots + f_n = 0
$$
Assume furthermore that for all $1 \le i_1 < i_2 < \cdots < i_s \le n$,
$$
f_{i_1} + f_{i_2} + \cdots + f_{i_s} = 0 \Longrightarrow 
\deg \gcd\{f_{i_1}, f_{i_2}, \ldots, f_{i_s}\} \le \deg f_n
$$
Let $d$ be the dimension of the vector space over $\C$ spanned by the
$f_i$s. Then
\begin{equation} \label{predav}
\max_{1 \le m \le n-1} \deg f_m - \deg f_n 
\le (d'-1)\left(\rg(f_1) + \rg(f_2) + \cdots + \rg(f_{n-1}) - \frac{d'}2\right)
\end{equation}
for all $d'$ between $d$ and $n-k'$ inclusive, where $k'$ is the number of constant $f_i$s with $i \le n-1$. Furthermore, equality in (\ref{predav}) is only possible if either $d' = 1$ (and therefore all $f_i$s are pairwise similar), or $\gcd\{f_1, f_2, \ldots, f_n\} = 1$.
\end{stelling}

\begin{bewijs}
Take $m' \le n-1$ such that $\max_{1 \le m \le n-1} \deg f_m = \deg f_{m'}$.
We reduce to the case that the conditions of theorem \ref{masons3} are satisfied.
If $f_n$ is constant, then the conditions of theorem \ref{masons3} are satisfied 
and hence we are done. So assume that $f_n$ is not constant. 

We first show that we may assume that 
$$
\deg f_n \le \min_{1 \le m \le n-1} \deg f_m
$$
We do this by removing all $f_i$s with $\deg f_i < \deg f_n$, adding them to $f_n$. This does not affect the $\gcd$-conditions on vanishing subsums, because for vanishing subsums with $f_n$, the $\gcd$-condition is automatically fulfilled. Notice further that $d$ can only decrease if we remove an $f_i$. So the smallest value of $d'$ will not increase. $n - k'$ will not change if we remove a constant $f_i$. But it will decrease $1$ if we remove a non-constant $f_i$. Therefore, for $d' = n - k'$, we cannot just take $d' = n - k'$ again in the reduced situation after removing a non-constant $f_i$. We take $d' = n - k' - 1$ instead. Since 
$$
-\frac{d'-1}2 < r(f_i) - \frac{d'}2
$$
this will not only work, but also makes \eqref{predav} strict. More generally, for all possible values of $d'$, \eqref{predav} will be strict if $d' \ge 2$ and there exists an $i$ such that $0 < \deg f_i < \deg f_n$.

We next show that we may assume that $\deg f_i = \deg f_n$ for at most one $i$ with $1 \le i \le {n-1}$. If $\deg f_i = \deg f_n$, then we cannot always remove $f_i$ in the above-described way, because we may have $\deg (f_n + f_i) < \deg f_n$. But if $\deg (f_n + f_{i'}) < \deg f_n$ as well, then $\deg (f_n + f_i + f_{i'}) = \deg f_n$. So we can remove $f_i$ and $f_{i'}$ simultaneously in that case.

Now we distinguish three cases. 
\begin{itemize}
 
\item
$n = 2$. \\
Then the left hand side of \eqref{predav} is zero and $d = 1$. The right hand side of \eqref{predav} is zero if $d' = 1$, and at least zero if $d' = 2$.
 
\item 
$n \ge 3$, and there is a minimal vanishing subsum of $f_1 + f_2 + \cdots + f_n = 0$ that contains both $f_{m'}$ and $f_n$ as summands. \\
Assume without loss of generality that $f_{m'} + f_{m'+1} + \cdots + f_n$ is a minimal vanishing subsum, and let $h := \gcd\{f_{m'},f_{m'+1},\ldots,f_n\}$. 

Notice that $\deg f_{m'} > \deg f_n$ and $d > 1$. If there is an $f_i$ with $m'+1 \le i \le n-1$ which is similar to $f_n$, then we can remove it in the above-described way, because $\deg f_n$ cannot decrease during the process. This way, we obtain that $f_i / h$ is not constant for $m' \le i \le n-1$. So
\begin{align*}
\lefteqn{\deg f_{m'} - \deg f_n} \\
&= \deg \frac{f_{m'}}{h} - \deg \frac{f_{n}}{h} \\
&\le (d' - 1) \left(\rg\left(\frac{f_{m'}}{h}\right) + 
      \rg\left(\frac{f_{m'+1}}{h}\right) + \cdots +
      \rg\left(\frac{f_{n-1}}{h}\right) - \frac{d'}2\right) \\
&\le (d' - 1) \left(\rg(f_{m'}) + \rg(f_{m'+1}) + \cdots + 
      \rg(f_{n-1}) - \rg(h) - \frac{d'}2\right) \\
&\le (d' - 1) \left(\rg(f_1) + \rg(f_2) + \cdots + \rg(f_{n-1}) - \frac{d'}2 - (m' - 1)\right) \\
&\le (d' - 1) \left(\rg(f_1) + \rg(f_2) + \cdots + \rg(f_{n-1}) - \frac{d'+m'-1}2\right)
\end{align*}
where $d'$ is at least the dimension of the vector space spanned by $f_{m'},\allowbreak f_{m'+1},\allowbreak \ldots,f_n$ and at most $n - m' + 1$. So \eqref{predav} follows. As $d' \ge d > 1$, $\deg h = 0$ and $m' = 1$ are required for equality to have a chance in \eqref{predav}.

\item 
$n \ge 3$ and there is no minimal vanishing subsum of $f_1 + f_2 + \cdots + f_n = 0$ that contains both $f_{m'}$ and $f_n$ as summands. \\
Assume without loss of generality that $f_1 + f_2 + \cdots + f_{m'}$ is a minimal vanishing subsum, and let $h := \gcd\{f_1,f_2,\ldots,f_{m'}\}$. Then $m' \le n-2$.

Notice again that $\deg f_{m'} > \deg f_n$ and $d > 1$. Since $\deg f_i = \deg f_n$ for at most one $i$ with $1 \le i \le m'$, we infer from $\deg h \le \deg f_n$ that $f_i / h$ is constant for at most one $i$ with $1 \le i \le m'$. By (\ref{kl}) in theorem \ref{masons3},
\begin{align*}
\lefteqn{\deg f_{m'} - \deg f_n} \\
&\le \deg \frac{f_{m'}}{h} \\
&\le (d' - 1) \left(\rg\left(\frac{f_1}{h}\right) + 
      \rg\left(\frac{f_2}{h}\right) + \cdots +
      \rg\left(\frac{f_{m'}}{h}\right) - \frac{d'}2\right) \\
&\le (d' - 1) \left(\rg(f_1) + \rg(f_2) + \cdots + \rg(f_{m'}) - r(h) - \frac{d'}2\right) \\
&\le (d' - 1) \left(\rg(f_1) + \rg(f_2) + \cdots + \rg(f_{n-1}) - \frac{d'}2 - (n-1-m')\right) \\
&\le (d' - 1) \left(\rg(f_1) + \rg(f_2) + \cdots + \rg(f_{n-1}) - \frac{d'+n-m'}2\right)
\end{align*}
where $d'$ is at least the dimension of the vector space spanned by $f_1,f_2,\ldots,\allowbreak f_{m'}$ and at most $m'$. So \eqref{predav} follows. As $d' \ge d > 1$, $m' = n$ is required for equality to have a chance in \eqref{predav}. But $m' \le n-2$, so equality is not possible in \eqref{predav}. \qedhere

\end{itemize}
\end{bewijs}

From the proof of theorem \ref{predavth}, we infer the following as well if not all $f_i$s are pairwise similar. If we remove all constant $f_i$s by adding them to $f_n$, and there is still a proper subsum of $f_1 + f_2 + \cdots + f_n$ which vanishes after that removal, then \eqref{predav} is strict before the removal of constant $f_i$s. The only case where this is not direct is when removing the constant $f_i$s results in that the $f_i$s become pairwise similar. But if the $f_i$s are pairwise similar and $n \ge 3$, then \eqref{predav} is strict for $d' \ge 2$. This is sufficient, because $d \ge 2$ before removing the constant $f_i$s.

Now substitute $f_i = g_i^{e_i}$ for all $i \le n-1$ and also
$f_n = -g_n = \sum_{i=1}^{n-1} g_i^{e_i}$, in (\ref{predav}). Then 
\begin{equation} \label{davenport}
\left(\sum_{i=1}^{n-1} \frac1{e_i} - \frac1{d'-1}\right)
\max_{1 \le m \le n-1} \deg g_m^{e_m}
\ge \frac{d'}2 - \frac1{d'-1} \deg \sum_{i=1}^{n-1} g_i^{e_i} 
\end{equation}
follows from (\ref{predav}) in a similar way as (\ref{fchy})
follows from (\ref{kl}) of theorem \ref{masons3}, see also \cite[(6)]{HY1}.

Indeed, applying (\ref{davenport}) on the sum $f^3 + ({\rm i}g)^2$ gives
$-\frac16 \deg(f^3) \ge 1 - \deg (f^3 - g^2)$ for $d' = 2$, 
which is equivalent to $\deg (f^3 - g^2) \ge \frac12 \deg f + 1$. 
For $d' = 3$, we get $\frac13 \deg f^3 \ge \frac32 - \frac12 \deg 
(f^3 - g^2)$, i.e.\@ $\deg(f^3 - g^2) \ge 3 - 2\deg f$, which is useless.

By replacing $n$ by $n+1$ in (\ref{davenport}), we obtain the main formula in the theorem below.

\begin{stelling}
Assume that $g_1 g_2 \cdots g_n$ is not constant, and
$$
g_1^{e_1} + g_2^{e_2} + \cdots + g_n^{e_n} \ne 0
$$
Assume furthermore that for all $1 \le i_1 < i_2 < \cdots < i_s \le n$,
$$
g_{i_1}^{e_{i_1}} + g_{i_2}^{e_{i_2}} + \cdots + g_{i_s}^{e_{i_s}} = 0 \Longrightarrow 
\deg \gcd\big\{g_{i_1}^{e_{i_1}}, g_{i_2}^{e_{i_2}}, \ldots, g_{i_s}^{e_{i_s}}\big\} \le \deg \sum_{i=1}^n g_i^{e_i}
$$
Let $d$ be the dimension of the vector space over $\C$ spanned by $g_1^{e_1}, g_2^{e_2}, \ldots, g_n^{e_n}$. Then
$$
\left(\sum_{i=1}^{n} \frac1{e_i} - \frac1{d'-1}\right)
\max_{1 \le m \le n} \deg g_m^{e_m}
\ge \frac{d'}2 - \frac1{d'-1} \deg \sum_{i=1}^n g_i^{e_i}
$$
for all $d'$ between $\max\{d,2\}$ and $n-k+1$ inclusive, where $k \ne n$ is the number of constant $g_i^{e_i}$s. Furthermore, equality cannot be reached in case $\gcd\{g_1, g_2, \ldots,\allowbreak g_n\} \ne 1$.
\end{stelling}

\begin{bewijs}
The result follows from theorem \ref{predavth}. If $\gcd\{g_1, g_2, \ldots,g_n\} \ne 1$, then $\sum_{i=1}^{n} g_i^{e_i}$ is not constant, and equality cannot be reached because $d' \ge 2$.
\end{bewijs}

In \cite{ZaD}, it is proved that for all even degrees of $f$, there are
univariate polynomials $f, g$ over $\C$ such that $\deg (f^3 - g^2) = 
\frac12 \deg f + 1$. Now assume $\deg (f^3 - g^2) = 
\frac12 \deg f + 1$. Then $\gcd\{f,g\} = 1$ and the Mason bound on
$-f^3 + g^2 + (f^3 - g^2) = 0$ gives us
$$
\deg f^3 \le \rg_1\Big(fg(f^3 - g^2)\Big) - 1 \le \deg \Big(fg(f^3 - g^2)\Big) - 1
$$
which is bound to be an equality. Furthermore, $fg(f^3 - g^2)$ is bound to be 
square-free. But any linear combination $\lambda f^3 + \mu g^2$ with 
$\lambda\mu \ne 0$ is bound to be square-free, since otherwise the inequality
$$
\deg f^3 \le \frac12 \Big(\rg_1(f^3) + \rg_1 (g^2) + \rg_1(f^3 - g^2) + 
             \rg_1(\lambda f^3 + \mu g^2) - 1\Big)
$$
would be violated. The above estimate is an instance of (\ref{dset}) in section
\ref{masondisc} below, since there exists a vanishing linear combination 
without zero coefficients of the arguments of $\rg_1$ on the right hand side.

\section{Some discussion on theorems \ref{masons3} and \ref{masons3a}}
\label{masondisc}

We describe now why the condition that all $f_i$s are relatively prime by pairs
is needed in \cite{BT,HY1,HY,ShSp}. They reduce to the case of maximal dimension
$d = n-1$ as follows. Assume that $f_n$ has the largest degree and say that 
$f_1, f_2, \ldots, f_d$ is a basis of the vector space over $\C$ spanned
by $f_1, f_2, \ldots, f_n$. Then 
$$
f_n = \lambda_1 f_1 + \lambda_2 f_2 + \cdots + \lambda_d f_d
$$
for some $\lambda_i \in \C$.
The greatest common divisor of the $f_i$s in the above sum is still the
same as in the original sum, but some $f_i$s might have a coefficient
$\lambda_i$ that is zero; say that $\lambda_1 \lambda_2 \cdots \lambda_{\rho} \ne 0$ 
and $\lambda_{\rho+1} = \lambda_{\rho+2} = \cdots = \lambda_d = 0$.
Then 
$$
\lambda_1 f_1 + \lambda_2 f_2 + \cdots + \lambda_{\rho} f_{\rho} + (-f_n) = 0
$$
is a vanishing sum of maximal dimension $\rho$. But the problem is
that the greatest common divisor of the the $f_i$s in the last sum
might be larger than that of the original sum.

But the above method does work when each set of $d$ $f_i$s generates
the whole vector space over $\C$ spanned by the $f_i$s, because that implies that 
$\rho = d$ above. So in this case 
one can get the estimates of theorems \ref{masons3} and \ref{masons3a}. 
But one can get estimates which are a factor $n-d$ smaller in this particular case, namely
\begin{equation} \label{dset}
\!\,\max_{1 \le m \le n} \deg f_m \le \frac1{n-d} \left( 
  \rg_{d-1}(f_1) + \rg_{d-1}(f_2) + \cdots + \rg_{d-1}(f_n)
  - \frac{d(d-1)}2 \right)\!\!\!\!\!
\end{equation}
and
\begin{equation} \label{dseta}
\max_{1 \le m \le n} \deg f_m \le \frac1{n-d} \left(
  \rg_{\frac{d(d-1)}{2}}(f_1 f_2 \cdots f_n) - \frac{d(d-1)}2 \right)
\end{equation}
combining techniques of \cite{HY} and the proof of \cite[Th.\@ 2]{Za}. 
We sketch the proof at the very end of this article.
  
In \cite[Th.\@ 2]{BB} it is shown that the coefficient $d'(d'-1)/2$ of
(\ref{kla}) in theorem \ref{masons3a} cannot be replaced by something less 
than $2n-5$, and the author conjectures that this
coefficient can indeed be improved to $2n-5$, i.e.
$$
\max_{1 \le m \le n} \deg f_m \le (2n-5) \big(\rg(f_1 f_2 \cdots f_n) - 1\big)
$$
I did not find similar considerations on (\ref{kl}) in theorem \ref{masons3}
in literature. So let us do something ourselves. The factor $(d'-1)$ in 
(\ref{kl}) cannot be improved, as is shown by the example
\begin{align*}
f_i &= \binom{n-2}{i-1} \left(x^{10^{100}}\right)^{i-1} 
        \quad (1 \le i < n) \\
f_n &= -\left(x^{10^{100}} + 1\right)^{n-2}
\end{align*}
The term $d'/2$ in (\ref{kl}) cannot be improved to $3d'/4$, 
as is shown by the example
\begin{align*}
f_i &= \lceil n/2 \rceil \binom{\lceil n/2 \rceil (\lfloor n/2 \rfloor + 1) - 2}{%
        \lceil n/2 \rceil i - 1} x^{\lceil n/2 \rceil i - 1} 
	\quad (i \le \lfloor n/2 \rfloor) \\
f_i &= - \zeta_{\lceil n/2 \rceil}^i 
        \left(x + \zeta_{\lceil n/2 \rceil}^i\right)^{\lceil n/2 \rceil (
	\lfloor n/2 \rfloor + 1) - 2} \quad (i > \lfloor n/2 \rfloor)
\end{align*}
for the case that none of the $f_i$s is constant, and by the example
\begin{align*}
f_i &= \lceil n/2 \rceil \binom{\lceil n/2 \rceil \lfloor n/2 \rfloor - 1}{%
        \lceil n/2 \rceil (i - 1)} x^{\lceil n/2 \rceil (i - 1)} 
	\quad (i \le \lfloor n/2 \rfloor) \\
f_i &= - \zeta_{\lceil n/2 \rceil}^{i} 
        \left(x + \zeta_{\lceil n/2 \rceil}^i\right)^{\lceil n/2 \rceil
	\lfloor n/2 \rfloor - 1} \quad (i > \lfloor n/2 \rfloor)
\end{align*}
for the case that $f_1$ is constant, but it might be possible to improve it to $3(d'-1)/\allowbreak 4$.

In section \ref{daven}, we have reduced (\ref{kl}) in theorem \ref{masons3} to
(\ref{klhy}) and (\ref{kla}) in theorem \ref{masons3a} to
(\ref{klahy}). Therefore it remains to prove
(\ref{klhy}) and (\ref{klahy}). But before we do that, we ask ourselves
the question whether (\ref{klhy}) and (\ref{klahy}) can be seen as instances
of one single, more general estimate. \cite{BrMa} has some valuable ideas in
that direction. Under the extra assumption that the $f_i$s are univariate
and $d = n-1$, (\ref{kla}) for $d' = d = n$ follows immediately from
\cite[Cor.\@ I]{BB}, and \cite[Cor.\@ II]{BB} implies
$$
\max_{1 \le m \le n} \deg f_m 
\le (n-2)\big(\rg(f_1) + \rg(f_2) + \cdots + \rg(f_n) + 1\big)
$$
but, since the $f_i$s are linearly independent, the number $k$ of constant 
$f_i$s is at most $1$. Since the number of empty $S_i$s in 
\cite[Cor.\@ II]{BB} equals $k$ as well, one can improve \cite[Cor.\@ II]{BB}
to
\begin{align} 
H(u_1,u_2,\ldots,u_n) &\le (n-2)\left(|S_1| + |S_2| + \cdots + |S_n| + k- 
\frac{n+1}2\right) - {} \nonumber \\ 
&\qquad \frac{(n-1)(n-2)}2 (2g - 2) \label{bbb}
\end{align}
and (\ref{kl}) in theorem \ref{masons3} for $d' = d = n$ follows.

The proof of (\ref{bbb}) is left as an exercise to the interested reader.
The general result that implies both \cite[Col.\@ I]{BB} and (the improved
version (\ref{bbb}) of) \cite[Col.\@ II]{BB} is \cite[Theorem A]{BB}. 

The rest of this article is organized as follows. In sections \ref{masonred}
to \ref{proofmain}, we prove (\ref{klhy}) of theorem \ref{masons3} and 
(\ref{klahy}) of theorem \ref{masons3a}. In section \ref{masonred}, we 
reduce to the univariate case. In section \ref{wrons}, we present the 
Wronskian, the key element in all generalized versions of Mason's theorem,
except \cite{Vo}. Section \ref{proofmain} consists of the actual proofs of
(\ref{klhy}) and (\ref{klahy}). At last, in section \ref{join}, we combine
(\ref{klhy}) and (\ref{klahy}) with ideas of \cite{BB}.

\section{Some reductions of the main theorem} \label{masonred}

By replacing the original sum by the minimal vanishing subsum containing $f_{m'}$
as a term, where $\deg f_{m'} = \max_{1 \le m \le n} \deg f_m$, we see that
in order to prove (\ref{klhy}) of theorem \ref{masons3}
and (\ref{klahy}) of theorem \ref{masons3a}, 
we can restrict ourselves to the case that $f_1 + f_2 + \cdots + f_n$ 
has no proper subsum that vanishes. 

We show now that we can restrict ourselves to the case that the $f_i$s are
univariate. More particular, a generic substitution $x_i = p_i y + q_i$ will
do the reduction. Assume that no proper
subsum of $f_1 + f_2 + \cdots + f_n$ vanishes and say 
that there are $l$ variables in the $f_i$s. Let $G$ be the set of
nonempty proper subsums
$$
f_{i_1} + f_{i_2} + \cdots + f_{i_s}
$$
and
$$
\bar{G} = \{ \bar g \mid g \in G \}
$$
where $\bar g$ is the largest degree homogeneous part of $g$ (i.e.\@ the sum
of all terms that have the same degree as $g$).	
Now pick a $p \in \C^l$ such that
$$
\bar{g}(p) \ne 0
$$
for all $\bar{g} \in \bar{G}$ (a $p$ that has coordinates that are 
transcendental over the field of coefficients of the $\bar{g}$s will do).

Assume without loss of generality that $p_1 \ne 0$ and define
$$
\hat{f}_i := f_i(p_1 x_1, x_2+p_2 x_1, \ldots, x_l+p_l x_1)
$$
for all $i$. Since $\gcd\{f_1, f_2, \ldots, f_n\} = 1$,
$\gcd\{\hat{f}_1, \hat{f}_2, \ldots, \hat{f}_n\} = 1$ as well. From Gauss's lemma, it follows that $\gcd\{f_1, f_2, \ldots, f_n\} = 1$ over $K(x_2,\ldots,x_l)[x_1]$ as well.
So if we apply the extended $\gcd$-theorem with respect to $x_1$,
we find $a_i \in \C(x_2,\ldots,x_l)[x_1]$ such that
$$
1 = a_1 \hat{f}_1 + a_2 \hat{f}_2 + \cdots + a_n \hat{f}_n
$$
For each $i$, write $a_i = \sum_{j=0}^{\infty} a_{i,j} x_1^j$ with 
$a_{i,j} \in \C(x_2,\ldots,x_l)$ and only finitely many $a_{i,j}$ nonzero.
Now put $q_1 := 0$ and take $(q_2, \ldots, q_l) \in \C^{k-1}$ such that
the denominators of the nonzero $a_{i,j}$s do not vanish on 
$(q_2, \ldots, q_l)$. Then
\begin{align}
1 &= a_1(q_2,\ldots,q_l)(x_1) \hat{f}_1(x_1,q_2,\ldots,q_l) + 
      {} \nonumber \\
  &\qquad a_2(q_2,\ldots,q_l)(x_1) \hat{f}_2(x_1,q_2,\ldots,q_l) + 
      \cdots + {} \nonumber \\    
  &\qquad a_n(q_2,\ldots,q_l)(x_1) \hat{f}_n(x_1,q_2,\ldots,q_l) \label{gcd1}
\end{align}
Put
$$
\tilde{f}_i := \hat{f}_i(y, q_2, \ldots, q_l) = f_i(q + yp) = 
               f_i(p_1 y + q_1,  p_2 y + q_2, \ldots, p_l y + q_l)
$$
for all $i$. From (\ref{gcd1}), it follows that
$\gcd\{\tilde{f}_1, \tilde{f}_2, \ldots, \tilde{f}_n\} = 1$.

Since $\rg_{\rho-1}(\tilde{f}_i) \le \rg_{\rho-1}(f_i)$ for all $i$ and
$\rg_{\sigma}(\tilde{f}_1 \tilde{f}_2 \cdots \tilde{f}_n) \le 
\rg_{\sigma}(f_1 f_2 \cdots f_n)$, it suffices to
show that $\deg \tilde{f}_i = \deg f_i$ for all $i$ and
no proper subsum of $\tilde{f}_1 + \tilde{f}_2 + \cdots + \tilde{f}_n = 0$
vanishes. We do so by proving that
for all proper subsets $I$ of $\{1,2,\ldots,n\}$:
$$
\deg \left(\sum_{i\in I} \tilde{f}_i\right) = 
   \deg \left(\sum_{i\in I} f_i\right)
$$
i.e.
$$
\deg g(q + yp) = \deg g
$$
for all $g \in G$. This is true, since the coefficient of $y^{\deg g}$
in $g(q + yp)$ is equal to $\bar{g}(p)$, which is nonzero by assumption.

\section{The Wronskian} \label{wrons}

Let $f_1, f_2, \ldots, f_n$ be polynomials in one and the same variable, 
say $y$. Then the {\em Wronskian determinant} 
of $f_1, f_2, \ldots, f_n$ is defined as
$$
W(f_1, f_2, \ldots, f_n) := \det \left( \begin{array}{cccc}
f_1    & f_2    & \cdots & f_n    \\
f_1'   & f_2'    & \cdots & f_n'   \\
\vdots & \vdots & \ddots & \vdots \\
f_1^{(n-1)} & f_2^{(n-1)} & \cdots & f_n^{(n-1)}
\end{array} \right)
$$
and the {\em Wronskian matrix} is the corresponding matrix 
on the right hand side.

Since differentiating is a linear operator, it
follows that $W(f_1, f_2, \ldots, f_n) = 0$ in case
\begin{equation} \label{lindepf}
\lambda_1 f_1 + \lambda_2 f_2 + \cdots + \lambda_n f_n = 0
\end{equation}
for some nonzero $\lambda \in \C^n$. Now a classical theorem tells
us that the reverse is true as well: if $f_1, f_2, \ldots, f_n$ are
linearly independent (i.e.\@ (\ref{lindepf}) implies $\lambda = 0$),
then $W(f_1, f_2, \ldots, f_n) \ne 0$. The example $f_1(x) = x^3$,
$f_2(x) = |x|^3$ shows us that the $f_i$s need to be polynomials.

Despite that the oldest known proof of this theorem by Frobenius is 
elementary, we give another proof, inspired by the proof of 
\cite[Lm.\@ 8]{Ye}. The reason for that will be given below. 

So let us assume that $f_1, f_2, \ldots, f_n$ are linearly independent. 
If there are two $f_i$s with the same degree, then we can subtract
a multiple of the first from the second to reduce the degree of the second,
since this operation does not affect the Wronskian determinant.
Progressing in this direction gives us that all $f_i$s have different 
degrees. Now order the $f_i$s by increasing degrees. This might only change
the sign of the Wronskian determinant.

The matrix
$$
\left( \begin{array}{cccc}
f_1^{(\deg f_1)} & f_2^{(\deg f_1)} & \cdots & f_n^{(\deg f_1)} \\
f_1^{(\deg f_2)} & f_2^{(\deg f_2)} & \cdots & f_n^{(\deg f_2)} \\
\vdots & \vdots & \ddots & \vdots \\
f_1^{(\deg f_n)} & f_2^{(\deg f_n)} & \cdots & f_n^{(\deg f_n)}
\end{array} \right)
$$
is upper triangular and does not have zeros on the diagonal.
Hence, its determinant does not vanish. Since it is a submatrix of 
$$
M := \left( \begin{array}{cccc}
f_1 & f_2 & \cdots & f_n \\
f_1' & f_2' & \cdots & f_n' \\
f_1^{(2)} & f_2^{(2)} & \cdots & f_n^{(2)} \\
\vdots & \vdots & \ddots & \vdots \\
f_1^{(\deg f_n)} & f_2^{(\deg f_n)} & \cdots & f_n^{(\deg f_n)}
\end{array} \right)
$$
this latter matrix has full rank $n$. 
Now we can make a square matrix $M'$ of full rank $n$ out of 
$M$ by throwing away redundant rows of $M$, i.e.\@ throwing away rows that 
are dependent of the rows above it. It suffices to prove that $M'$
is the Wronskian matrix, i.e.\@
$$
M' = \left( \begin{array}{cccc}
f_1 & f_2 & \cdots & f_n \\
f_1' & f_2' & \cdots & f_n' \\
f_1^{(2)} & f_2^{(2)} & \cdots & f_n^{(2)} \\
\vdots & \vdots & \ddots & \vdots \\
f_1^{(n-1)} & f_2^{(n-1)} & \cdots & f_n^{(n-1)}
\end{array} \right)
$$
Write $f^{(i)}$ for the vector
$$
(f_1^{(i)}, f_2^{(i)}, \cdots, f_n^{(i)})
$$
and $f = f^{(0)}$ and $f' = f^{(1)}$. Assume that the $m$-th row of $M'$ is 
$(f^{(m-1)})\tp$, but the $(m+1)$-th row of $M'$ is not $(f^{(m)})\tp$, say it is 
$(f^{(j)})\tp$ with $j > m$. Then $(f^{(j-1)})\tp$ is in the space generated by 
the first $m$ rows of $M'$, i.e.
\begin{equation} \label{fj1}
f^{(j-1)} = a_0 f + a_1 f' + a_2 f^{(2)} + \cdots + a_{m-1} f^{(m-1)}
\end{equation}
where the $a_i$ are rational functions, i.e.\@ quotients of polynomials,
for all $i$. Differentiating (\ref{fj1}) gives
$$
f^{(j)} = (a_0' f + a_0 f') + (a_1' f' + a_1 f^{(2)}) + \cdots + 
	  (a_{m-1}' f^{(m-1)} + a_{m-1} f^{(m)})
$$
Since each of the $2m$ terms on the right hand side is contained in the 
space generated by the first $m$ rows of $M'$, $f^{(j)}$ is contained in this
space as well. Contradiction, so the $m$-th row of $M'$ is $(f^{(m-1)})\tp$ for
all $m$.

In \cite[Lemma 6, pp.\@ 15-16]{Sch}, a generalization of the Wronskian theorem
for more variables is formulated. The operators 
$\frac{\partial^i}{\partial y^i}$ are in fact replaced by operators
$\Delta_i$, each of which is a product of partial derivatives.
The number of partial derivatives that $\Delta_i$ decomposes into, 
multiple appearances counted by their frequency, is called the {\em order}
$o(\Delta_i)$ of $\Delta_i$. 

The usual Wronskian determinant is replaced by
\begin{equation} \label{genwrons}
W_{\Delta}(f_1, f_2, \ldots, f_n) := \det \left( \begin{array}{cccc}
\Delta_1 f_1 & \Delta_1 f_2 & \cdots & \Delta_1 f_n \\
\Delta_2 f_1 & \Delta_2 f_2 & \cdots & \Delta_2 f_n \\
\vdots & \vdots & \ddots & \vdots \\
\Delta_n f_1 & \Delta_n f_2 & \cdots & \Delta_n f_n
\end{array} \right)
\end{equation}
and the author T. Schneider of \cite{Sch} 
proves that if $f_1, f_2, \ldots, f_n$ are
linearly independent, then $W_{\Delta}(f_1, f_2, \ldots, f_n) \ne 0$ for
certain operators $\Delta_i$ of order $i-1$ at most. In particular,
$\Delta_1$ is the identity operator, and the first row looks the same
as in the case of one variable.

Unlike the above proof of the classical Wronskian theorem,
the proof of this theorem by Frobenius cannot be generalized to
more indeterminates. The way Schneider proves his multivariate result is
by reducing to the univariate Wronskian theorem. But his theorem
does not show that there are $\Delta_i$s of all orders $0,1,2,\ldots,\rho$, 
where $\rho$ is the maximum order of the $\Delta_i$s, unlike a 
straightforward generalization of the above proof of the classical 
Wronskian theorem to more indeterminates. Neither does his methods give tools to
prove that
\begin{equation} \label{gcdwrons}
W_{\Delta}(hf_1, hf_2, \ldots, hf_n) = h^n W_{\Delta}(f_1, f_2, \ldots, f_n)
\end{equation}
(\ref{gcdwrons}) can be found in \cite[Lm.\@ 2.1]{HY}. But this lemma
is somewhat different to both our methods and \cite[Lemma 6, pp.\@ 15-16]{Sch}, 
since the Wronskian determinant might be zero.
 
Take for instance $f = (1, xy, x^2y^2)$. Notice that
$$
W_{1,\frac{\partial}{\partial x},\frac{\partial^2}{\partial x^2}}
(1, xy, x^2y^2) = \det \left( \begin{array}{ccc} 
1 & xy & x^2y^2 \\
0 & y & 2xy^2 \\
0 & 0 & 2y^2 
\end{array} \right) = 2y^3
$$
and this is also a generalized Wronskian
one can get by the multivariate variant of the above method, since
$\frac{\partial}{\partial y} x^i y^i = i x^i y^{i-1} = 
\frac{x\partial}{y\partial x} x^i y^i$. The above Wronskian matrix is however 
{\em not} of the form of \cite[Lm.\@ 2.1]{HY} and \cite[Lm.\@ 8]{Ye}, 
because $\frac{\partial}{\partial y} f$ is {\em not} linearly dependent over $\C$ of
its rows. The Wronskian matrix of both lemma's must be that of
$$
W_{1,\frac{\partial}{\partial x},\frac{\partial}{\partial y}}(1, xy, x^2y^2) = 0
$$
instead.

In the proofs of theorems \ref{masons3} and \ref{masons3a}, we shall
employ a special generalized Wronskian, one without an identity operator:

\begin{lemma} \label{wronslem}
Let $\tilde{f}_1, \tilde{f}_2, \ldots, \tilde{f}_n$ be polynomials over $\C$ in the variables
$y,z_1,\allowbreak z_2,\allowbreak \ldots,z_l$ which are linear in $z_1,z_2,\ldots,z_l$.
Assume that $\tilde{f}_1, \tilde{f}_2, \ldots, \tilde{f}_n$ are linearly independent over $\C$. 
Then there exists a $\Delta = (\Delta_1, \Delta_2, \ldots, \Delta_n)$ with
$$
W_{\Delta}(\tilde{f}_1,\tilde{f}_2,\ldots,\tilde{f}_n) \ne 0
$$
such that for each $i$, either
$$
\Delta_i = \frac{\partial}{\partial z_j}
$$
for some $j$, or (if $i \ge 2$)
$$
\Delta_i = \frac{\partial}{\partial y} \Delta_{i-1}
$$
\end{lemma}

\begin{bewijs}
Choose $j$ such that $\lambda_{j,n} \ne 0$. Say that
$\lambda_{j,1} = \cdots = \lambda_{j,m} = 0$ and
$\lambda_{j,m+1} \cdots \lambda_{j,n} \ne 0$.
We distinguish three cases:
\begin{itemize}

\item $\frac{\partial}{\partial z_j} \tilde{f}_{m+1}, 
\ldots, \frac{\partial}{\partial z_j} \tilde{f}_n$ are linearly dependent. \\
After possibly interchanging $\tilde{f}_{m+1}$ with one of $\tilde{f}_{m+2}, \ldots, \tilde{f}_{n-1}$, we can write
$$
\frac{\partial}{\partial z_j} \tilde{f}_{m+1} = 
\mu_{m+2} \frac{\partial}{\partial z_j} \tilde{f}_{m+2} + \cdots + 
\mu_n \frac{\partial}{\partial z_j} \tilde{f}_n
$$
Replace $\tilde{f}_{m+1}$ by $\tilde{f}_{m+1} - (\mu_{m+2} \tilde{f}_{m+2} + \cdots + 
\mu_n \tilde{f}_n)$ and apply induction on $n-m$.

\item $\frac{\partial}{\partial z_j} \tilde{f}_{m+1}, 
\ldots, \frac{\partial}{\partial z_j} \tilde{f}_n$ are linearly independent
and $m = 0$. \\
Then the result follows by applying the Wronskian theorem (in one variable)
on $\frac{\partial}{\partial z_j} \tilde{f}_1, \frac{\partial}{\partial z_j} \tilde{f}_2,
\ldots, \frac{\partial}{\partial z_j} \tilde{f}_n$. The operators
are $\Delta_i = \frac{\partial^i}{\partial y^{i-1} \partial z_j}$.

\item $\frac{\partial}{\partial z_j} \tilde{f}_{m+1}, 
\ldots, \frac{\partial}{\partial z_j} \tilde{f}_n$ are linearly independent
and $m \ge 1$. \\
From the above case, it follows that $W_D(\tilde{f}_{m+1},\ldots,\tilde{f}_n) \ne 0$,
where $D_i := \frac{\partial^i}{\partial y^{i-1} \partial z_j}$. By induction on $n$, we have $W_{\Delta} (\tilde{f}_1, \tilde{f}_2, \ldots, \tilde{f}_m) \ne 0$.
Now extend $\Delta$ by defining $\Delta_{m+i} = D_i$ for all $i \ge 1$.
Since $\frac{\partial}{\partial z_j} \tilde{f}_{i} = 0$ for all $i \le m$, it follows that
$$
W_{\Delta} (\tilde{f}_1, \tilde{f}_2, \ldots, \tilde{f}_n) = 
W_{\Delta} (\tilde{f}_1, \ldots, \tilde{f}_m) \cdot W_D (\tilde{f}_{m+1}, \ldots, \tilde{f}_n) \ne 0
$$
and $\Delta$ remains of the desired form. \qedhere

\end{itemize}
\end{bewijs}

Notice that the above lemma can be generalized to more variables
as well.

\section{Proof of the main theorem} \label{proofmain}

From the reductions in sections \ref{daven} and \ref{masonred}, 
it follows that in order to prove theorems \ref{masons3} and \ref{masons3a}, 
it suffices to prove the following:

\begin{stelling} \label{masontilde}
Let ${f}_1, {f}_2, \ldots, {f}_n$ be nonzero polynomials
over $\C$ in the variable $y$, not all constant, such that 
$\gcd\{{f}_1, {f}_2, \ldots, {f}_n\} = 1$ and
$$
{f}_1 + {f}_2 + \cdots + {f}_n = 0
$$
Let $d$ be the dimension of the vector space over $\C$ spanned by
the ${f}_i$s and assume furthermore that no proper subsum of
${f}_1 + {f}_2 + \cdots + {f}_n$ vanishes. Then
$$
\max_{1 \le m \le n} \deg {f}_m \le 
   \rg_{\rho-1}({f}_1) + \rg_{\rho-1}({f}_2) + \cdots + 
   \rg_{\rho-1}({f}_n) - \frac{\rho(\rho-1)}2
$$
for some $\rho$ with $2 \le \rho \le d$, and
$$
\max_{1 \le m \le n} \deg {f}_m \le 
   \rg_{\sigma}({f}_1{f}_2 \cdots {f}_n) - \sigma 
$$
for some $\sigma$ with $1 \le \sigma \le d(d-1)/2$.
\end{stelling}

Assume without loss of generality that ${f}_1, {f}_2, \ldots, 
{f}_d$ is a basis of the vector space over $\C$ spanned by the
${f}_i$s. For each $j > d$, there exists unique $\lambda_{j,i}$
such that
\begin{equation} \label{minrel}
{f}_j = \sum_{i=1}^d \lambda_{j,i} {f}_i
\end{equation}

In order to get rid of all linear relations between the ${f}_i$s
except the sum relation, we define
$$
\tilde{f}_i := \left(\sum_{j=d+1}^n \lambda_{j,i} z_j\right) \cdot {f}_i
$$
for all $i \le d$, and
$$
\tilde{f}_i := -z_i \cdot {f}_i
$$
for all $i > d$. It follows from (\ref{minrel}) that
\begin{align*} 
\sum_{i=1}^n \tilde{f}_i &= \sum_{i=1}^d \sum_{j=d+1}^n \lambda_{j,i} z_j {f}_i 
                     - \sum_{j=d+1}^n z_j {f}_j \\
                 &= \sum_{j=d+1}^n z_j \left(\sum_{i=1}^d \lambda_{j,i} 
		     {f}_i - {f}_j \right) \\
		 &= 0
\end{align*}

Furthermore, it follows from (\ref{minrel}) that
$$
\sum_{i=1}^d \left(1 + \sum_{j=d+1}^n \lambda_{j,i}\right) {f}_i
= \sum_{i=1}^d {f}_i + 
  \sum_{j=d+1}^n \sum_{i=1}^d \lambda_{j,i} {f}_i
= \sum_{i=1}^n {f}_i = 0
$$
whence 
\begin{equation} \label{lambdasum}
\sum_{j=d+1}^n \lambda_{j,i} = -1 \quad (1 \le i \le d)
\end{equation}
for ${f}_1, {f}_2, \ldots, {f}_d$ are linearly 
independent.

\begin{lemma} \label{mu1}
$\mu_1 \tilde{f}_1 + \mu_2 \tilde{f}_2 + \cdots + \mu_n \tilde{f}_n = 0$ implies
$\mu_1 = \mu_2 = \cdots = \mu_n$.
\end{lemma}

\begin{bewijs}
Let $G$ be the graph with vertices
$\{1,2,\ldots,n\}$ and connect two vertices $j,i$ by an edge if 
$\lambda_{j,i} \ne 0$. Notice that $G$ is a bipartite graph
between $\{1,2,\ldots,d\}$ and $\{d+1,\ldots,n\}$.
We first show that $G$ is connected. Assume the opposite.
Say that $G$ does not have an edge between $\{1,\ldots,d',d+1,\ldots,n'\}$ and 
$\{\mbox{$d'+1$},\allowbreak \ldots,d,n'+1,\ldots,n\}$, where either $d' < d$ or $n' < n$. 
Then $\lambda_{j,i} = 0$ for all $j > n'$ and $i \le d'$, whence by (\ref{lambdasum})
\begin{equation} \label{lambdasum1}
\sum_{j=d+1}^{n'} \lambda_{j,i} = -1 
\end{equation}
for all $i \le d'$. On the other hand, $\lambda_{j,i} = 0$ for all $j \le n'$ 
and $i > d'$, whence
\begin{equation} \label{lambdasum2}
\sum_{j=d+1}^{n'} \lambda_{j,i} = 0
\end{equation}
for all $i > d'$.

Substituting $z_j = 1$ for all $j \le n'$ and $z_j = 0$ for all $j > n'$
in $\sum_{i=1}^n \tilde{f}_i$, it follows from (\ref{lambdasum1}) and (\ref{lambdasum2}) 
that we obtain
$$
\sum_{i=1}^{d} \left( \sum_{j=d+1}^{n'} \lambda_{j,i} \right) 
{f}_i - \sum_{j=d+1}^{n'} {f}_j
= -\sum_{i=1}^{d'} {f}_i - \sum_{j=d+1}^{n'} {f}_j
$$
which is zero, since $\sum_{i=1}^n \tilde{f}_i$ is zero. Since no proper subsum of 
$\sum_{i=1}^n {f}_i$ vanishes, we have $d' = d$ and $n' = n$. Contradiction, so
$G$ is connected.

Now assume $\mu_1 \tilde{f}_1 + \mu_2 \tilde{f}_2 + \cdots + \mu_n \tilde{f}_n = 0$. Pick a $j > d$.
Substituting
$z_j = 1$ and $z_m = 0$ for all $m \ne j$ in $\sum_{i=1}^n \mu_i \tilde{f}_i$ gives us
$$
\sum_{i=1}^d \mu_i \lambda_{j,i} {f}_i - \mu_j {f}_j = 0
$$
but on account of (\ref{minrel}), also
$$
\sum_{i=1}^d \mu_j \lambda_{j,i} {f}_i - \mu_j {f}_j = 0
$$
so by subtraction
$$
\sum_{i=1}^d (\mu_i-\mu_j) \lambda_{j,i} {f}_i = 0
$$
Since ${f}_1, {f}_2, \cdots, {f}_d$ are linearly independent
over $\C$, $(\mu_i-\mu_j) \lambda_{j,i} = 0$ for all $i \le d$. So
\begin{equation} \label{muk}
\lambda_{j,i} \ne 0 ~ \Longrightarrow ~ \mu_i = \mu_j
\end{equation}
Since $G$ is connected, the desired result follows.
\end{bewijs}

From lemma \ref{mu1}, it follows that $\tilde{f}_1, \tilde{f}_2, \ldots, \tilde{f}_{n-1}$ 
are linearly independent, whence we can
apply lemma \ref{wronslem} to get
$$
W_{\Delta} (\tilde{f}_1, \tilde{f}_2, \ldots, \tilde{f}_{n-1}) \ne 0
$$
where $\Delta = (\Delta_1, \Delta_2, \ldots, \Delta_{n-1})$ satisfies
the properties of lemma \ref{wronslem}. Since 
$\tilde{f}_1 + \tilde{f}_2 + \cdots + \tilde{f}_n = 0$, we have
\begin{align}
W_{\Delta} (\tilde{f}_1, \tilde{f}_2, \ldots, \tilde{f}_{n-1}) 
&= (-1)^{n-i} W_{\Delta} (\tilde{f}_1, \ldots, \tilde{f}_{i-1}, \tilde{f}_{i+1}, \ldots, \tilde{f}_n) 
    \nonumber \\
&= (-1)^{n-1} W_{\Delta} (\tilde{f}_2, \ldots, \tilde{f}_{n-1}, \tilde{f}_n) \label{rhowrons}
\end{align}
Let $\rho$ be the maximum among the orders $o(\Delta_1), o(\Delta_2), 
\ldots, o(\Delta_{n-1})$, i.e.\@ the maximum number of partial derivatives
which any $\Delta_m$ may decomposes into. Put
$$
\sigma := \sum_{i=1}^{n-1} (o(\Delta_i) - 1)
$$

Since $\frac{\partial}{\partial z_j} \tilde{f}_n = 0$ for all $j \ne n$, and the right hand side of (\ref{rhowrons}) does not vanish, $\frac{\partial}{\partial z_n} \in \{\Delta_1, \Delta_2, \allowbreak \ldots, \Delta_{n-1}\}$. A similar argument on the left hand side of (\ref{rhowrons}) gives $\frac{\partial}{\partial z_m} \in \{\Delta_1, \Delta_2, \ldots, \Delta_{n-1}\}$ for $m$ with $d < m < n$. So $n-d$ of the $n-1$ $\Delta_i$s have order $1$. From the conditions of theorem \ref{masontilde}, we infer that $d \ge 2$. It follows from lemma \ref{wronslem} that 
$$
2 \le \rho \le d \qquad \mbox{and} \qquad 
1 \le \frac{\rho(\rho-1)}2 \le \sigma \le \frac{d(d-1)}2
$$

\begin{lemma} \label{divlem}
$$
{f}_1 {f}_2 \cdots {f}_n \,\Big|\,
   \rd_{\rho-1}({f}_1) \rd_{\rho-1}({f}_2) \cdots 
   \rd_{\rho-1}({f}_n) \cdot W_{\Delta} (\tilde{f}_1, \tilde{f}_2, \ldots, \tilde{f}_{n-1})
$$
and
$$
{f}_1 {f}_2 \cdots {f}_n \,\Big|\,
   \rd_{\sigma}({f}_1{f}_2 \cdots {f}_n) \cdot 
   W_{\Delta} (\tilde{f}_1, \tilde{f}_2, \ldots, \tilde{f}_{n-1})
$$
\end{lemma}

\begin{bewijs}
It suffices to prove that irreducible polynomials $g$ over $\C$ in the
variable $y$ divide the right hand side at least as often as the left hand 
side. So let $g \in \C[y]$ be irreducible. Since 
$\gcd\{{f}_1, {f}_2, \ldots, {f}_n\} = 1$, 
one of the ${f}_i$s is not divisible by $g$, say that
$g \nmid {f}_1$. It follows from (\ref{rhowrons}) that it
suffices to show that $g$ divides ${f}_2 \cdots {f}_n$ 
at most as often as 
$$
\rd_{\rho-1}({f}_1) \rd_{\rho-1}({f}_2) \cdots 
   \rd_{\rho-1}({f}_n) \cdot W_{\Delta} (\tilde{f}_2, \ldots, \tilde{f}_{n-1}, \tilde{f}_n)
$$
and
$$
\rd_{\sigma}({f}_1{f}_2 \cdots {f}_n) \cdot 
   W_{\Delta} (\tilde{f}_2, \ldots, \tilde{f}_{n-1}, \tilde{f}_n)
$$

Now pick any term of the determinant expression
$W_{\Delta} (\tilde{f}_2, \ldots, \tilde{f}_{n-1}, \tilde{f}_n)$. After permuting $\tilde{f}_2, \ldots, \tilde{f}_n$,
the term at hand becomes
$$
\Delta_1 \tilde{f}_2 \cdot \Delta_2 \tilde{f}_3 \cdot \cdots \cdot \Delta_{n-1} \tilde{f}_n
$$
Now if $g$ divides ${f}_i$ exactly $l$ times and hence also $\tilde{f}_i$ 
exactly $l$ times, then $g$ divides $\Delta_{i-1} \tilde{f}_i$ at 
least $l - \rho$ times, since partial derivatives kill at most one instance
of a factor $g$ in their argument. But one of the partial derivatives is
a $\frac{\partial}{\partial z_j}$ which does not kill any instance of $g$, so
$g$ divides $\Delta_{i-1} \tilde{f}_i$ at least $l - (\rho - 1)$ times.

The factor $\rd_{\rho-1}({f}_i)$
compensates the decrease of $\rho-1$ factors $g$, so $g$ divides
$\rd_{\rho-1}({f}_i) \Delta_{i-1} {f}_i$ at least as often 
as it divides ${f}_i$, and the first inequality of this lemma follows.
The second inequality follows from the fact that the $\Delta_i$s together
have $\sigma$ partial derivatives of the form $\frac{\partial}{\partial y}$
that might kill instances of $g$.
\end{bewijs}

\begin{lemma} \label{deglem}
$$
\deg W_{\Delta} (\tilde{f}_1, \tilde{f}_2, \ldots, \tilde{f}_{n-1}) 
   \le \deg ({f}_1{f}_2 \cdots {f}_{n-1}) - \sigma
$$
\end{lemma}

\begin{bewijs}
The idea is that a partial derivative decreases the degree by one.
Consider a term on the left hand side of the above formula. After 
reordering the $\tilde{f}_i$s, this term becomes
$$
\Delta_1 \tilde{f}_1 \cdot \Delta_2 \tilde{f}_2 \cdot \cdots \cdot \Delta_{n-1} \tilde{f}_{n-1} 
$$
Since $o(\Delta_i) \ge 1$ for all $i$, the degree of this term is at most
$\deg (\tilde{f}_1\tilde{f}_2 \cdots \tilde{f}_{n-1}) - (n-1) = 
\deg ({f}_1{f}_2 \cdots {f}_{n-1})$. But there are also
$\Delta_i$s of orders larger than one, which are responsible for the
term $\sigma$.
\end{bewijs}

\begin{bewijs}[Proof of theorem \ref{masontilde}]
Assume without loss of generality that ${f}_n$ has the largest degree 
among the ${f}_i$s. From lemmas \ref{divlem} and \ref{deglem}, 
it follows that
$$
\sum_{i=1}^n \deg {f}_i 
   \le \rg_{\rho-1}({f}_1) + \rg_{\rho-1}({f}_2) + \cdots +
       \rg_{\rho-1}({f}_n) + 
       \deg ({f}_1{f}_2 \cdots {f}_{n-1}) - \sigma
$$
whence
$$
\deg {f}_n 
   \le \rg_{\rho-1}({f}_1) + \rg_{\rho-1}({f}_2) + \cdots +
       \rg_{\rho-1}({f}_n) - \frac{\rho(\rho-1)}2
$$
which is the first inequality of theorem \ref{masontilde}. The second
inequality follows similarly.
\end{bewijs}
 
\section{Joining theorems \ref{masons3} and \ref{masons3a}} \label{join}

The general result that implies both \cite[Col.\@ I]{BB} and (the improved
version (\ref{bbb}) of) \cite[Col.\@ II]{BB} is \cite[Theorem A]{BB}, which
we will describe now for the polynomial case. 
For irreducible polynomials $p$, let $m_p$ denote the number of $f_i$s that
is {\em not} divisible by $p$. Then \cite[Theorem A]{BB} implies
\begin{equation} \label{thA}
\max_{1 \le m \le n} \deg f_m \le  -\binom{n-1}{2} + 
\sum_p \left( \binom{n-1}{2} - \binom{m_p-1}{2} \right)
\end{equation}
where $\sum_p$ ranges over all irreducible polynomials $p$. 
It follows from (\ref{thA}) that
$$ 
\max_{1 \le m \le n} \deg f_m \le -\binom{n-1}{2} + 
\sum_{p \nmid f_1 \cdots f_n} \left(\binom{n-1}{2} - \binom{n-1}{2}\right)
+ \sum_{p \mid f_1 \cdots f_n} \binom{n-1}{2}
$$
which is exactly the case $d'= n-1$ of the univariate case of 
(\ref{kla}) in theorem \ref{masons3a}.

In order to get a similar result on (\ref{thA}) and (\ref{kl}) in theorem 
\ref{masons3}, we first need some preparations. Assume
\begin{equation} \label{nodiv}
f_i \nmid f_{i+1}
\end{equation}
The reason for (\ref{nodiv}) is that there exists an irreducible $p$ that divides
$f_i$ more times than it divides $f_{i+1}$, say that $p$ divides
$f_i$ $l+j$ times and $f_{i+1}$ $l$ times. Now replace
$f_i$ by $f_i p^{j}$ and $f_{i+1}$ by $f_{i+1} p^{-j}$. Then (\ref{nodiv})
might still be the case, but the divisibility by $p$ is not the reason any 
more. Furthermore, for any power $q$ of an irreducible polynomial,
$q$ divides as many $f_i$s as before. If we proceed in this direction,
we finally arrive at

\begin{propositie} \label{hprop}
There exist $h_1, h_2, \ldots, h_n$ such that
\begin{enumerate}

\item $h_1 \mid h_2 \mid \cdots \mid h_n$,

\item For any power $q$ of an irreducible polynomial,
      $q$ divides as many $h_i$s as it divides $f_i$s.
      
\end{enumerate}
\end{propositie}

Notice that $h_1 = \gcd\{f_1, f_2, \ldots, f_n\} = 1$. More generally,
$h_i$ is the greatest common divisor over all subsets
$\{j_1, j_2, \ldots, j_i\}$ of $\{1,2,\ldots,n\}$ of
$\lcm\{f_{j_1}, \allowbreak f_{j_2}, \allowbreak \ldots, \allowbreak f_{j_i}\}$.

Since $m_p$ is also the number of $h_i$s that is not divisible
by $p$,
$$
\binom{m_p - 1}{2} = \sum_{i=2}^{m_p} (i-2) =
\sum_{\genfrac{}{}{0pt}{}{2 \le i \le n}{p \nmid h_i}} (i-2)
$$
whence
$$
\binom{n-1}{2} - \binom{m_p - 1}{2} = \sum_{i=m_p+1}^n (i-2) =
\sum_{\genfrac{}{}{0pt}{}{2 \le i \le n}{p \mid h_i}} (i-2)
$$
Summing this over all $p$, it follows from (\ref{thA}) that
\begin{equation} \label{thaa}
\max_{1 \le m \le n} \deg f_m \le \sum_{i=2}^n (i-2) \rg(h_i) - \binom{n-1}{2}
\end{equation}
which implies the case $d'= n-1$ of the univariate case of 
(\ref{kl}) in theorem \ref{masons3}, for
$$
\sum_{i=1}^n \rg(h_i) = \sum_{i=1}^n \rg(f_i)
$$

By $\rg(h_1) = 0$ and $\rg(h_i) \le \rg(h_n)$, the case $d'= n-1$ of the univariate 
case of (\ref{kla}) in theorem \ref{masons3a} follows from (\ref{thaa}) as 
well. (\ref{thaa}) can be improved to 
$$
\max_{1 \le m \le n} \deg f_m \le \sum_{i=3}^n \rg_{i-2}(h_i) - \binom{n-1}{2}
$$
which implies (\ref{klhy}) in theorem \ref{masons3} for $\rho = n-1$ and
(\ref{klahy}) in theorem \ref{masons3a} for $\sigma = (n-1)(n-2)/2$,
since $\rg_i(a) \rg_j(b) \le \rg_{i+j}(ab)$. The general multivariate result 
that includes both theorems \ref{masons3} and \ref{masons3a} is as follows.

\begin{stelling}
Let $f_1, f_2, \ldots, f_n$ be (possibly multivariate) nonzero polynomials
over $\C$ in the variable $y$, not all constant, such that 
$\gcd\{f_1, f_2, \ldots, f_n\} = 1$ and
$$
f_1 + f_2 + \cdots + f_n = 0
$$
Let $d$ be the dimension of the vector space over $\C$ spanned by
the $f_i$s and assume furthermore that no proper subsum of
$f_1 + f_2 + \cdots + f_n$ vanishes. 

Take $h_1, h_2, \ldots, h_n$ as in proposition \ref{hprop}. Then
$$
\max_{1 \le m \le n} \deg f_m \le \left( \rg_{(o_1)-1}(h_2) + \rg_{(o_2)-1}(h_3) + \cdots + \rg_{(o_{n-1})-1}(h_n) \right) - \sigma
$$
where
$$
o_1 \le o_2 \le \cdots \le o_{n-1}
$$
are the orders of the $\Delta_i$s. 
\end{stelling}

\begin{proof}
The proof is similar to that of theorem \ref{masontilde}, except for lemma \ref{divlem}. Instead, we need
$$
f_1 f_2 \cdots f_n \mid \rd_{o_1-1}(h_2)\rd_{o_2-1}(h_3)\cdots\rd_{o_{n-1}-1}(h_n) W_{\Delta} (\tilde{f}_1, \tilde{f}_2, \ldots, \tilde{f}_{n-1})
$$
In a term of the Wronskian determinant expression, each operator $\Delta_i = \frac{\partial}{\partial y}^j \frac{\partial}{\partial z_k}$ kills at most $j$ factors $p$ of its operand, but only to the extend that this operand has factors $p$. The number of factors $p$ which is killed by operators $\Delta_i$ in a Wronskian determinant term is bounded from above by the number of factors $p$ which are killed by operators $\frac{\partial}{\partial y}$ in
$$
\bigg(\frac{\partial}{\partial y}^{o_1-1} h_2\bigg)\bigg(\frac{\partial}{\partial y}^{o_2-1} h_3\bigg)\cdots\bigg(\frac{\partial}{\partial y}^{o_{n-1}-1} h_n\bigg)
$$
From this, the proof can be obtained.
\end{proof}

At last we sketch the proof of (\ref{dset}) and (\ref{dseta}). Assume that
each set of $d$ $f_i$s forms a basis of the space generated by all
$f_i$s, and order the $f_i$s by increasing degree. As indicated in
section \ref{masondisc}, we do not need to multiply the $f_i$s by linear 
forms in order to get rid of unwanted linear dependences. 
Similar to (\ref{rhowrons}), one can prove that all sequences of $d$ 
$f_i$s have the same Wronskian determinant $W_{\Delta}(f_1, f_2, \ldots, f_d)$ 
up to a nonzero constant in $\C$. Since each set of $d$
$f_i$s generates the whole space, the greatest common divisor of such a 
set is $1$, whence there can only be $d-1$ $f_i$s at most that are divisible
by a given irreducible polynomial $p$. So $h_1 = h_2 = \cdots = h_{n-d+1} = 1$,
$f_1 f_2 \cdots f_n = h_{n-d+2} h_{n-d+3} \cdots h_n$, and
\begin{equation}
f_1 f_2 \cdots f_n \mid \rd_1(h_{n-d+2}) \rd_2(h_{n-d+3}) \cdots \rd_{d-1}(h_n)
       W_{\Delta}(f_1, f_2, \ldots, f_d) \label{dseth}
\end{equation}
because focusing on one irreducible divisor $p$, one can replace
$f_1, f_2, \ldots, f_d$ on the right hand side of (\ref{dseth}) by the $d$ 
$f_i$s of maximum divisibility by $p$. Next, since each set of $d$
$f_i$s has an element of maximum degree and the $f_i$s are ordered 
by increasing degree, we infer that
$$
\deg f_d = \deg f_{d+1} = \cdots = \deg f_n
$$
From that, we can deduce the factor $1/(n-d)$ in (\ref{dset}) and (\ref{dseta}).

\end{document}